\documentclass[numbers,sort&compress]{IntechOpen-Book}%%OPTIONS are numbers, authoryear, sort&compress,sectref

\graphicspath{{Artworks/}}
\begin{document}

\Mainmatter

\begin{frontmatter}

\chapter{A Review of the EnKF for Parameter Estimation}
% Enter author names EXACTLY as you would like them to appear in the final manuscript; only use pattern: first name, last name. Do NOT use pattern: last name, first name. Patronyms fall into the first name category.
\author{Neil K. Chada}
%\author{Simon Weissmann}
%\author{Author 3}

\makechaptertitle

\chaptermark{Ensemble Kalman Inversion}

\begin{abstract} % abstract = max 200 words, unstructured format, single paragraph
The ensemble Kalman filter is a well-known and celebrated data assimilation algorithm. It is of particular relevance as it used for high-dimensional problems,
by updating an ensemble of particles through a sample mean and covariance matrices. In this chapter we present a relatively recent topic which is the application
of the EnKF to inverse problems, known as ensemble Kalman Inversion (EKI). EKI is used for parameter estimation, which can be viewed as a black-box optimizer
for PDE-constrained inverse problems. We present in this chapter a review of the discussed methodology, while presenting emerging and new areas of research, where 
numerical experiments are provided on numerous interesting models arising in geosciences and numerical weather prediction. 
\end{abstract}

\begin{keywords} % use a minimum of 5 kwrds, separate them with a comma
Ensemble Kalman filter, Kalman filter, inverse problems, parameter estimation, data assimilation, optimization
\end{keywords}

\end{frontmatter}

\section{Introduction} % first section MUST be titled Introduction, and feature introductory text; do NOT change this title

Inverse problems \cite{KS04, AT87, AMS10} are a class of mathematical problems which have gained significant attention of recent. Simply put, inverse problems are concerned with the recovery of some parameter
of interest from noisy unstructured data. Mathematically we can express an inverse problem as the recovery of $u\in \mathcal{X}$ from noisy measurements of data $y \in \mathcal{Y}$,
expressed as
\begin{equation}
\label{eq:ip}
y = \mathcal{G}(u) + \eta,
\end{equation} 
where $\mathcal{G}:\mathcal{X} \rightarrow \mathcal{Y}$ is the forward operator, and $\eta \sim \mathcal{N}(0,\Gamma)$ is some form of additive Gaussian noise. Inverse problems are of 
high interest due to the amount of relevant problems that arise in wide variety of applications, most notably geophysical sciences, medical imaging and numerical weather prediction \cite{ORL08,ACL86,MW06}. The classical
approach to solving inverse problems, which is the theme of this chapter, is to construct a least-squares functional, and the solution is represented as a minimizer of some functional of the form
\begin{equation}
\label{eq:lsf}
u^*:= \arg \min_{u \in \mathcal{X}} \frac{1}{2}\|y - \mathcal{G}(u)\|^2_{\Gamma} + \lambda R(u),
\end{equation}
where $\lambda>0$ is a regularization parameter and $R(u)$ is some regularization term, usually required to prevent the overfitting of the data. Traditional methods 
for solving \eqref{eq:ip} include optimization schemes such as the Gauss--Newton method, or Levenburg--Marquardt method which require derivative information of $\mathcal{G}$,
which can prove costly and cumbersome. Therefore a motivation for solving inverse problems is to provide gradient-free optimizers which can reduce this computational
burden, while attaining a good level of accuracy. The methodology that we motivate, which alleviates these issues, is that of ensemble Kalman inversion (EKI). EKI can be
viewed as the application of the ensemble Kalman filter (EnKF) to inverse problems, which is a natural way to solve inverse problems given the connections between data assimilation
and inverse problems. The EnKF is a Monte-Carlo version of the celebrated Kalman filter, which is more favorable in high-dimensions. It operates by updating an ensemble of particles
through sample mean and covariances. In particular we will take the viewpoint of EKI which acts as PDE-constrained derivative-free optimizer. Therefore EKI can be viewed as a black-box
solver where no derivative information is required. Since this method was proposed for inverse problems, it has seen wide applications to various engineering-based applications,
as well as developments related to both theory and methodology. In this chapter we discuss some of these keys concepts and insights, while briefly mentioning particular directions
with EKI. 

The general outline of these chapter is as follows. In Section \ref{sec:EKI_back} we provide the necessary background material, which covers the basics of EKI 
with some intuition and motivation We will discuss the algorithm in both the usual discrete-time setting, but also the continuous-time setting. This will lead onto Section 
\ref{sec:reg} where we discuss one recent direction which is that of regularization theory, and its application to EKI. Furthermore we will also discuss how EKI
can be extended to the notion of sampling in statistics within Section \ref{sec:samp}. Other, less-developed, directions are provided in Section \ref{sec:other}. 
 Numerical experiments are provided in basic settings in Section \ref{sec:num} on a number of basic differential equations, before providing some future remarks and a conclusion in Section \ref{sec:conc}.

\section{EKI: Background material}
\label{sec:EKI_back}

%The body is where the author explains experiments, presents and interprets data of one's research. Authors are free to decide how the main body will be structured. However, you are required to have \textbf{at least one heading}. Please ensure that either British or American English is used consistently in your chapter.

%The text throughout the manuscript will be \textbf{left-aligned (or ragged-right)} in the final version of the chapter. This is not a typesetting error. This cannot be changed on an individual basis, i.e. IntechOpen will not accept requests for custom text alignment. All chapters, in all publications, will have the same layout and formatting.

In this section we provide the background material related to the understanding and intiution of EKI. This will begin with a discussion on the ensemble Kalman filter,
and how it connections with EKI. We will then present EKI in its vanilla form, which is a discrete-time optimizer, before discussing its connections
with various existing methods. Finally we will extend the original formulation to the setting of continuous-time where we aim to provide a gradient
flow structure of the resulting equations.

\subsection{Kalman filtering}

The ensemble Kalman filter (EnKF), is a popular methodology based on the celebrated Kalman filter (KF), which was originally developed my 
Rudolph Kalman in the 1960s \cite{REK60,RSB65}. The Kalman filters initial aim was to solve a recursive estimation problem from dynamics processes
and systems. Specifically the KF aims to merge data with model, or signal, dynamics where both equations have the form
\begin{align}
\label{eq:signal_da}
u_{n+1}&=\Psi(u_{n}) + \xi_{n}, \  \ \ \ \ \ \  \ \ \{\xi_{n}\}_{n \in \mathbb{Z}^{+}} \sim \mathcal{N}(0,\Sigma), \\ 
\label{eq:data_da}
y_{n+1}&=H(u_{n+1}) + \eta_{n+1}, \ \ \  \ \{\eta_{n+1}\}_{n \in \mathbb{Z}^{+}} \sim \mathcal{N}(0,\Gamma).
\end{align}
Here $\{u_{n}\}_{n \in \mathbb{Z}^{+}}$ is our signal which is updated through a forward operator $\Psi:\mathbb{R}^m \rightarrow \mathbb{R}^m$, which when combined with noise,  provides the update $u_{n+1}$. Our data is denoted as $y_{n+1}$ which is produced by sending our updated signal through the operator $H:\mathbb{R}^m \rightarrow \mathbb{R}^{\overline{m}}$, where $\overline{m} >m$, which is known as observational operator. Our initial conditions for the system are given as $u_0 \sim \mathcal{N}(m_0,\mathcal{C}_0)$. This area of recursive estimation, in this setup, became to be known as data assimilation \cite{BC09,LSZ15}.
\\ \bigskip
In particular in the linear and Gaussian setting, where the dynamics and noise are Gaussian, the KF updates state 
using the first two moments, which we know are the mean and covariance. Assume that the state-space dimension is $d \in {R}^+$,
then the cost of the KF has complexity $\mathcal{O}(d^2)$. For high-dimensional examples this can be an issue, therefore an algorithm 
that was developed to alleviate this is the EnKF, a Monte Carlo version, proposed by Evensen \cite{GE09,GE03}.

The EnKF operates by replacing the true covariance by a sample covariance and mean and updates an ensemble of particles $u^{(j)}_n$, with $1\leq j \leq J$ particles,
using these moments combined with information from the data.
The EnKF can be split into a two-step procedure, which is the prediction step

\begin{align*}
\hat{u}^{(j)}_{n+1} &= \Psi(u^{(j)}_{n}) + \xi^{(j)}_{n}, \quad
\hat{m}_{n+1} = \frac{1}{J} \sum^{J}_{j=1} u^{(j)}_{n+1}, \\
\hat{{C}}_{n+1} &=\frac{1}{J-1} \sum^{J}_{j=1} (u^{(j)}_{n+1} - \hat{m}_{n+1} )(u^{(j)}_{n+1} - \hat{m}_{n+1} )^{T},
\end{align*}
and update step
\begin{align*}
K_{n+1} &= \hat{{C}}_{n+1}H^T(H\hat{{C}}_{n+1}H^T + \Gamma), \\
{u}^{(j)}_{n+1} &=(I - K_{j+1}H)\hat{u}^{(j)}_{n+1} + K_{n+1}y^{(j)}_{n+1}, \\
y^{(j)}_{n+1} &= y_{n+1} + \eta^{(j)}_{n+1},
\end{align*}
where $K_{n+1}$ represents the Kalman gain matrix and $\xi^{(j)}_{n}$ and $\eta^{(j)}_{n+1}$ are i.i.d. Gaussian noise.
 In the EnKF context our prediction step defines a sample mean and covariance from our signal. From this in the analysis step we define our Kalman gain 
through our sample covariance, which updates our signal, which is given by ${u}^{(j)}_{n+1}$. This is aided by aiming to minimize the discrepancy of the data 
$y^{(j)}_{n+1}$ and the quantity $H(u)$.  To better understand this discrepancy, there is an alternative approach of looking at the EnKF is through a variational approach, 
where we consider the follow cost function
\begin{equation}
\label{eq:cf}
I_n(u) := \frac{1}{2}|y^{(j)}_{n+1} - H(u)|^2_{\Gamma} + \frac{1}{2}|u - \hat{u}_{n+1}^{(j)}|^2_{\hat{{C}}_{n+1}},
\end{equation}
for which we aim to minimise, which is defined as the updated mean
\begin{equation}
\label{eq:va}
\hat{m}_{n+1} = \arg \min_{u}I_n(u).
\end{equation}
This minimization procedure relies on the updated covariance $\hat{C}_{n+1}$ which is dependent entirely on $\hat{u}^{(j)}$. As described in the {prediction step} and {update step} of filtering, a mapping is presented between distributions. As we related the distributions in the filtering setting, for each step, we can do so similarly for the EnKF, i.e. 
\begin{equation*}
\{u^{(j)}_n\}_{j=1}^J \mapsto \{u^{(j)}_{n+1}\}_{j=1}^J, \ \ \ \ \ \ \ \{{u}^{(j)}_{n+1}\}_{j=1}^J \mapsto \{\hat{u}^{(j)}_{n+1}\}_{j=1}^J.
\end{equation*}
With the EnKF, compared to KF, the computational complexity associated with it is $\mathcal{O}(Jd)$, where
one usually assumes $J<d$, therefore implying the reduction in cost.

\subsection{EnKF applied to inverse problems}

Since the formulation of the EnKF, there has been a huge interest from practitioners in various applicable disciplines.
Most notably this has been within numerical weather prediction, geophysical sciences and signal processing related
to state estimation. In this chapter our focus is on the application of the EnKF to inverse problems, namely to solve
\eqref{eq:ip}. We now introduce this application which is known as ensemble Kalman inversion (EKI), which was
introduced by Iglesias et al., motivated from Li et al, \cite{ILS13} as a derivative-free optimizer for PDE-constrained inverse problems.

\begin{figure}[!h]\centering
	\FIG{\includegraphics[width=0.7\textwidth,height=12pc]{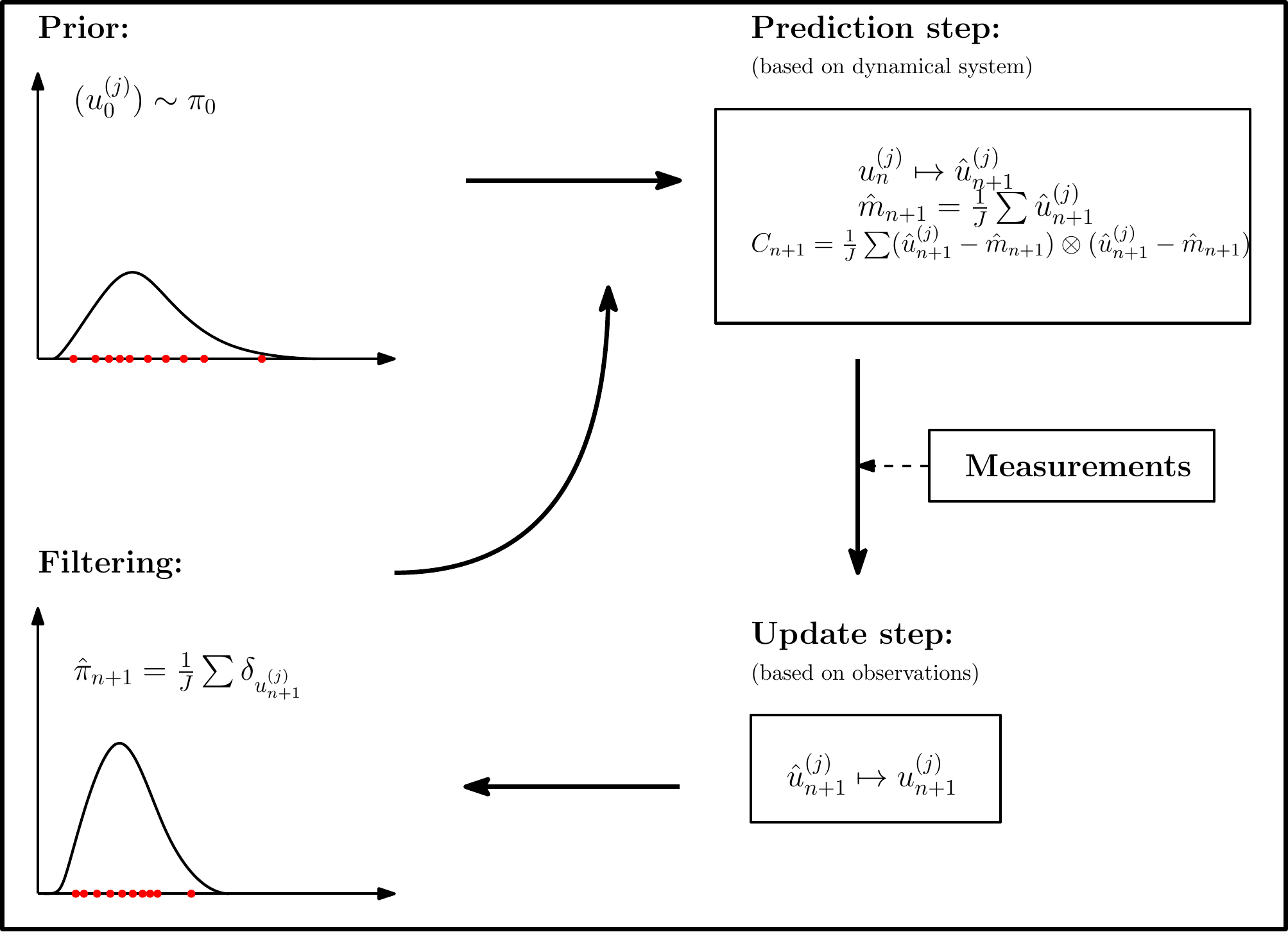}}
	{\caption{Dynamics of the ensemble Kalman filter, split into the prediction and update steps.}}
	\label{fig:enkf}
\end{figure}

As with the EnKF, we are concerned with updating an ensemble of particles, for which now we modify notation with $n$ now denoting
the iteration count. Given an initial ensemble $\{u^{(j)}_0\}$, our aim is to learn a true underlying
unknown $u^{\dagger}$. To do so, as done with the EnKF, we first define our sample mean and covariance matrices 
\begin{align*}
&\bar{u}^{(j)}_n = \frac{1}{J} \sum^{J}_{j=1} u^{(j)}_n, \quad \bar{u}^{(j)}_n = \frac{1}{J} \sum^{J}_{j=1} {G}(u^{(j)}_n), \\
C^{uu}_n = \frac{1}{J-1}\sum^J_{j=1}(u^{(j)}_n &- \bar{u})(u^{(j)}_n - \bar{u})^T, \quad C^{up}_n = \frac{1}{J-1}\sum(u^{(j)}_n - \bar{u})({G}(u^{(j)}_n) - \bar{{G}})^T.
\end{align*}
which we can through the update equation
\begin{align}
\label{eq:update}
u^{(j)}_{n+1} &= u^{(j)}_n + hC^{up}(hC^{pp} + \Gamma)^{-1}(y^{(j)}_n - \mathcal{G}(u^{(j)}_n)), \\
\label{eq:data}
y^{(j)}_n &= y + \eta^{(j)}_n,
\end{align}
where $y$ represents our true data and $h>0$ denotes a step size related to the level of discretization.
\textbf{Figure \ref{fig:enkf}} provides a pictorial description of the EnKF, which has been described above.

The update equation of EKI \eqref{eq:update} is of interest as it coincides with the update formula for
Tikhonov regularization for linear statistical inverse problems. Namely if we consider $R(u) = \frac{1}{2}\|u\|^2_{C_0}$,
then the update formula, in the linear $\mathcal{G}(\cdot) = \mathcal{G}$ and Gaussian setting is given as 
\begin{equation}
\label{eq:tp}
{u}_{\textrm{TP}} =  \bar{{u}} + C\mathcal{G}^*(\mathcal{G}C\mathcal{G}^*+\Gamma)^{-1}({y} - \mathcal{G} \bar{{u}}), 
\end{equation}
where $\mathcal{G}^*$ denotes the derivative of the operator $\mathcal{G}$. This connection is of relevance and was
discussed in \cite{ILS13}, where it was shown that taking the limit as $J \rightarrow \infty$, it was shown that $u \rightarrow u_{TP}$.
This is of interest as the minimizing the regularized functional \eqref{eq:tp} is equivalent to the following maximization
procedure in statistics
$$
u:= \arg \max_{u \in \mathcal{X}}  \mathbb{P}(u|y).
$$
known as the MAP formulation, where   $\mathbb{P}(u|y) = \mathbb{P}(y|u)\mathbb{P}(u)$ denotes the posterior distribution
. This connection is discussed in \cite{LPS89}. Therefore this provides some insight into EKI and its connection
with other known existing methodologies in inverse problems. An important entity to discuss is a property that EKI
inherits, which is the \textit{subspace property}. It is given by the following lemma.

\begin{lemma}
Let $\mathcal{A}$ be the linear span of the initial ensemble $\{{u}_0^{(j)}\}_{j=1}^J$, then we that
$\{\textcolor{black}{u}^{(j)}_n\}_{j=1}^J \in \mathcal{A}$ for all $n\in\mathbb{N}$. 
\end{lemma}

The essence of the subspace property states that the updated ensemble of particles is spanned by
the initial ensemble. This is important, because it provides a justification on the performance, whether
the initial ensemble is a good choice or not. Therefore it can act as an advantage or a disadvantage.

\subsection{Continuous-time formulation}\label{ssec:cont_time}

The original representation of EKI, as shown in \eqref{eq:update} , is a discrete-time iterative scheme similar to other optimization methods.
However it is of interest to understand EKI in a continuos-time setting, which was considered by Schillings et al. \cite{SS17,BCW19}. This is primarily for two reasons; (i) firstly that one
can understand more easily how the dynamics of \eqref{eq:update}-\eqref{eq:data} behaves, and secondly (ii) it provides new
numerical schemes for EKI, which is specific in the continuous-time setting. In order to derive such equations, as usual we 
require to take the step-size to zero, i.e. $h \rightarrow 0$. Once we do this, we have the following set of stochastic differential
equations
\begin{equation*}
\frac{du^{(j)}}{dt} = C^{uw}(u)\Gamma^{-1}\big(y - \mathcal{G}(u^{(j)})\big) + C^{uw}(u)\sqrt{\Gamma^{-1}}\frac{dW^{(j)}}{dt},
\end{equation*}
with $W^{(j)}$ denoting independent cylindrical  Brownian motions. By substituting the form of the covariance operator, we see
\begin{equation}
\label{eq:gen1}
\frac{du^{(j)}}{dt} = \frac{1}{J}\sum^{J}_{k=1} \Big\langle \mathcal{G}(u^{(k)}) - \bar{\mathcal{G}},y- \mathcal{G}(u^{(j)}) + \sqrt{\Gamma}\frac{dW^{(j)}}{dt} \Big\rangle_{\Gamma}(u^{(k)} - \bar{u}).
\end{equation}

For this we take our forward operator $\mathcal{G}(\cdot) = A\cdot$ to be bounded and linear. Using this notion and by substituting our linear operator $A$ in \eqref{eq:gen1} we have the following diffusion limit
\begin{equation}
\label{eq:linear}
\frac{du^{(j)}}{dt} = \frac{1}{J}\sum^{J}_{k=1} \Big\langle A(u^{(k)}- \bar{u}),y- Au^{(j)} \Big\rangle_{\Gamma}(u^{(k)} - \bar{u}).
\end{equation} 
By defining the empirical covariance operator 
\begin{equation*}
C(u) = \frac{1}{J-1}\sum^{J}_{k=1}(u^{(k)} - \bar{u}) \otimes (u^{(k)} - \bar{u}),	
\end{equation*}
and taking $\Gamma=0$ we can express \eqref{eq:linear} as 
\begin{align}
\label{eq:gf-nh}
\frac{du^{(j)}}{dt} &= -C(u) D_{u} \Phi(u^{(j)};y), \\
\Phi(u;y)&= \frac{1}{2} \| \Gamma ^{-1/2}(y - Au)\|^2. \nonumber
\end{align}
Thus we note that each particle performs a preconditioned gradient descent for $\Phi(\cdot;y)$ where all the gradient descents are preconditioned through the covariance $C(u)$. Since our covariance operator $C(u)$ is semi-positive definite we have that
\begin{equation*}
\frac{d}{dt}\Phi(u(t);y) = \frac{d}{dt}\frac{1}{2} \| \Gamma ^{-1/2}(y - Au)\|^2 \leq 0.
\end{equation*}
In the context of EKI this is of interest as it is a first result providing some indication of the dynamics, which was not achievable through the discrete-time
update formula \eqref{eq:update}. Indeed given the gradient flow structure, we are able to see that the EKI abides by a usual optimization function, with the dynamics
following the direction of the negative gradient, or in other-words towards to minimizer of $\Phi$. Since the continuous-time formulation was derived, there 
has been different works deriving further analysis, most notably with recent success on the nonlinear setting, and other well-known results. This can be found 
in \cite{BSW21}.

\section{Regularization}
\label{sec:reg}

\begin{figure}[!h]\centering
	\FIG{\includegraphics[width=0.5\textwidth,height=9pc]{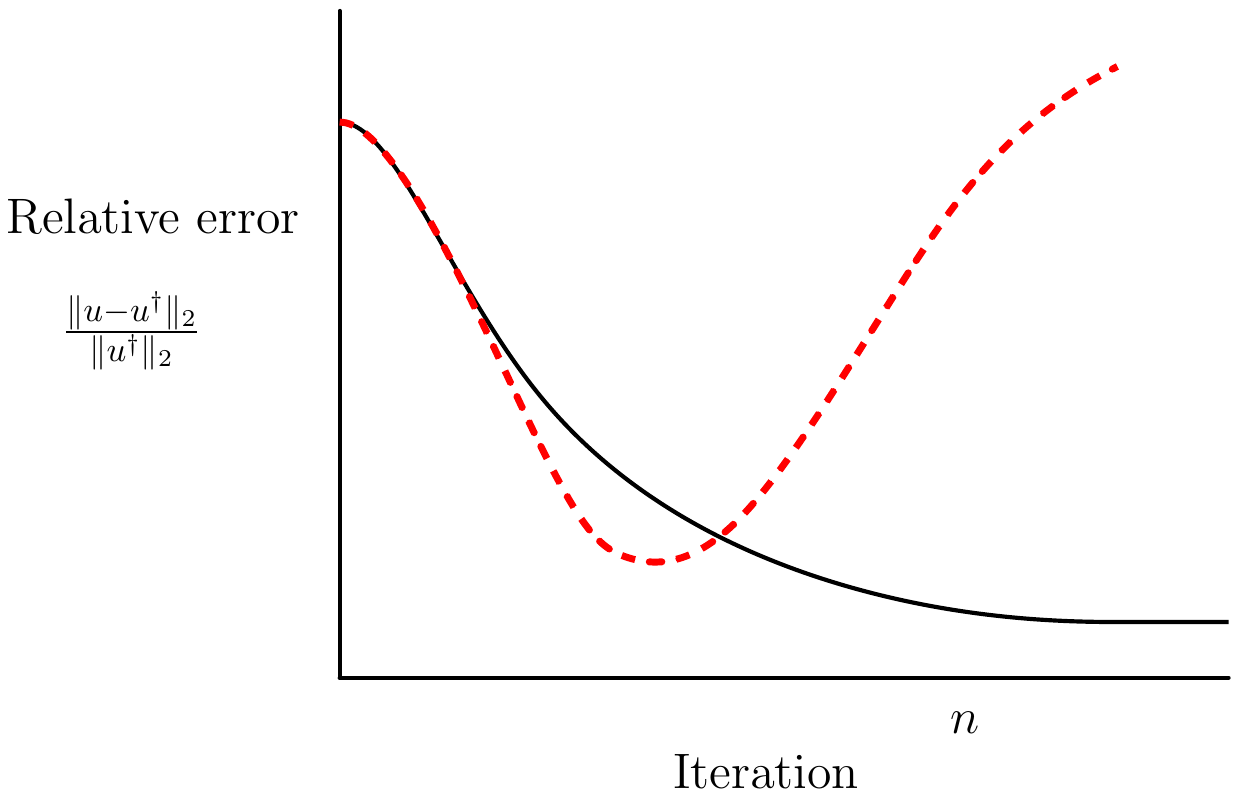}}
	{\caption{The figure presents two simulations of EKI as the iterations increase. The black curve represents what we aim to achieve, however in
	certain situations the data is commonly overfitted. Therefore this can cause a divergence in the relative error, as shown by the dashed red curve.}}
	\label{fig:reg}
\end{figure}

In this section we discuss the role of regularization in EKI. We will begin with an introduction into iterative regularization schemes, that have been used
before discussing Tikhonov regularization, $L_p$ and particular adaptive choices.
\\ \bigskip
As briefly discussed regularization is an important tool in optimization, and inverse problems
aimed at preventing the over-fitting, or influence, of the data. We refer the reader to various pieces of literature that give a concise overview on this \cite{BB18,EHN96}.
The over-fitting of data can cause issues in inverse problems, such as the divergence of the error, therefore careful consideration is needed to prevent this.
A cartoon representation of this is given in \textbf{Figure \ref{fig:reg}}.

To initiate this chapter, there are two main forms of regularization one can apply for inverse problems. The first is related to \textit{iterative regularization},
where the regularization is included within the iterative scheme. This can be included directly such as the form
\begin{equation}
\label{eq:up2}
u^{(j)}_{n+1} = u^{(j)}_n + hC^{up}(hC^{pp} + \alpha_n\Gamma)^{-1}(y^{(j)}_n - \mathcal{G}(u^{(j)}_n)), 
\end{equation}
or in the presence of a discrepancy principle of the form
$$
\|\Gamma^{-1}(y-\bar{u}^{(j)}_n ) \|^2 \leq \vartheta \eta, \qquad \vartheta \in (0,1),
$$
which controls the error between the updated ensemble and the true unknown. The discrepancy principle acts as a stopping rule if the error becomes big,
and the the modified update formula contains a sequence of numbers $\{\alpha_n\}_{n \in \mathbb{N}}$ aimed at also preventing the overfitting of the data.
This sequence is chosen in such a way that is related to a discrepancy principle. Specifically for EKI this has been considered in numerous work by 
Iglesias et al. \cite{MAI16,IY20}.

However more recent work has considered regularization through the least-squares functional (LSF) \eqref{eq:lsf}. For EKI the first known form to consider this, is
Tikhonov regularization which has the penalty form of $R(u) = \frac{1}{2}\|u\|^2_{C_0}$. This form of regularization is a natural choice, as it very well-known
and understood but can view viewed as a Gaussian form of regularization, which smoothes the problems. In the context of EKI this makes sense, as commonly
one assumes Gaussian dynamics. The work of Chada et al. \cite{CST19} first developed this extension, which was done by modifying \eqref{eq:ip} to the following

\begin{subequations}
\begin{align*}
\textcolor{black}{y}&=G(\textcolor{black}{u}) + \eta_1, \\ 
\textcolor{black}{u}&= \eta_2, 
\end{align*}
\end{subequations}
where $\eta_1 \sim \mathcal{N}(0,\Gamma), \eta_2 \sim \mathcal{N}(0,\lambda^{-1}C_0).$
\\ \bigskip
Now we introduce $z,\eta$ and the mapping $\mathcal{F}: \mathcal{X} \times \mathcal{X} \mapsto
\mathcal{Y} \times \mathcal{X}$ as follows:
\[
\textcolor{black}{z}=\begin{bmatrix}
\textcolor{black}{y}\\
0
\end{bmatrix},\quad
F(u)=\begin{bmatrix}
{G}(\textcolor{black}{u})\\
\textcolor{black}{u}
\end{bmatrix},
\quad
\eta=\begin{bmatrix}
\eta_1\\
\eta_2
\end{bmatrix},
\]
and
\[
\eta \sim N(0,\Sigma), \quad
\Sigma =
\begin{bmatrix}
\Gamma & 0\\
0 &\lambda^{-1}C_0
\end{bmatrix}.
\]
Therefore our inverse problem is now reformulated at 
\begin{equation*}
  \textcolor{black}{z} = \mathcal{F}(\textcolor{black}{u}) + \eta.
\end{equation*}
now from this we can modify EKI to include the above setup, for which we refer to it
as Tikhonov ensemble Kalman inversion (TEKI), which takes the following form
$$
u^{(j)}_{n+1} = u^{(j)}_n + hB^{up}(hB^{pp} + \Gamma)^{-1}(z^{(j)}_n - \mathcal{F}(u^{(j)}_n)), 
$$
where we have now modified covariance matrices $B^{up}, B^{pp}$.
From this inclusion, the authors of \cite{CST19} were able to show that analytically, the subspace property 
still holds, while other such results as observability and controllability and the ensemble collapse.
More importantly through the numerical simulations, it was shown that one can prevent the over-fitting
phenomenon. 

Since this work a number of useful extensions have been considered, such as its understanding in the continuous-case,
as well as the new variants in the discrete-time setting \cite{CT21}. Two recent developments on this have been firstly
on the extension to $L_p$ regularization \cite{YL21,SSW20}, which is to motivate reconstructing edges or lines, where the LSF is modified
to
$$
\Phi(u;y) := \frac{1}{2}\|y - \mathcal{G}(u)\|^2_{\Gamma} + \lambda \|u\|_{p}, \quad p \geq 1.
$$
Finally another direction is related to producing adaptive strategies for TEKI. Adaptive regularization schemes are of importance,
as choosing a correct choice of the regularization parameter $\lambda>0$ can have a big impact on the reconstruction. Therefore
thinking adaptively allows one to evolve the parameter over the iteration count, now denoted as $\lambda_n$. The work of Weissmann et al. \cite{WCS22}
provides these developments in an adaptive fashion.

\section{Ensemble Kalman sampling}
\label{sec:samp}
%\swtd{Check if all abbreviations are introduced}

Although the EKI has been introduced through the application of the EnKF to inverse problems and hence sequential sampling method,  the trending viewpoint of EKI lies in optimization. So far, we have seen its motivation from the gradient flow structure in the continuous-time formulation in Section~\ref{ssec:cont_time} and the representation as SDE. For applying EKI as a consistent sampling method, we would instead of taking the limit $t\to\infty$ rather consider the limit $t\to1$. For linear forward models EKI is consistent with the posterior distribution, however,  it is known to be not consistent with the Bayesian perspective in the nonlinear setting \cite{ESS2015}. 

Building up on this fact, the motivation behind the ensemble Kalman sampler \cite{GFLA20} is to modify the time-dynamical system of EKI in  a way such that the limiting distribution for $t\to\infty$ corresponds to the posterior distribution.  We will start the discussion with an introductory example.

%%%%%%%%%%%%%
\begin{example}
Let $\pi_\ast$ be a pdf of the form $\pi_\ast(x)\propto \exp(-\Phi(u))$ with $\Phi(u) = \frac12\|y-G(u)\|_\Gamma^2 + \|u\|_{C}^2$, i.e.~$\pi_\ast$ corresponds to the posterior pdf under Gaussian prior assumption $\pi_0=\mathcal N(0,C)$. We consider the Langevin diffusion given by
\begin{equation}\label{eq:Langevin}
du_t = \nabla_u\log\pi_\ast(u_t)\, dt + \sqrt{2}dW_t, \quad u_0\sim\pi_0,
\end{equation}
where $(W_t)_{t\ge0}$ denotes a Brownian motion in $\mathcal X = \mathbb R^{n_u}$.  The evolution of the distribution $\rho_t$ of the state $u_t$ can then be described through the Fokker--Planck equation
\begin{equation}\label{eq:FP}
\partial \rho_t = \nabla\cdot (\rho_t \nabla\log\pi_\ast) + \Delta \rho_t, \quad \rho_0 = \pi_0,
\end{equation}
where under certain assumptions on $\Phi$ the underlying Markov process $(u_t)_{t \geq 0}$ is ergodic and its unique invariant distribution is given by $\pi_{\ast}$ \cite{P2014}.  Taking the Fokker--Planck equation \eqref{eq:FP} into account the convergence to equilibrium can be described through the Kullback--Leibler (KL) divergence $\mathrm{KL}= \int_{\mathcal{X}} q_1(x)\log(\frac{q_1(x)}{q_2(x)})\,{\mathrm{d}}x$ \cite{KL1951}.
Assuming a $\log$-Sobolev inequality (e.g.~satisfied for $\log$-concave $\pi_\ast$), it follows that
$$
\mathrm{KL}(\rho_t\mid\pi_\ast)\leq  \exp(-\lambda t) \mathrm{KL}(\rho_0\mid\pi_\ast)
$$
for some $\lambda>0$ \cite{MV1999}.
\end{example}

\subsection{Interacting Langevin sampler}

The interacting Langevin sampler has been introduced, motivated by the preconditioned gradient descent method, as interacting particle system represented by the coupled system of SDEs
\[du_t^{(j)} = C(u_t)\nabla_u\log\pi_\ast(u_t^{(j)})\, dt + \sqrt{2C(u_t)}dW_t, \quad j=1,\dots,J,\]
initialized through an i.i.d.~sample $u_0^{(j)}\sim\pi_0$. The idea of preconditioning with $C(u_t)$ instead of a fixed preconditioning matrix $C\in\mathbb{R}^{n_u\times n_u}$ is motivated through the corresponding mean-field limit. In the large particle limit, the corresponding SDE is given as
\[du_t = C(\rho_t)\nabla_u\log\pi_\ast(u_t)\, dt + \sqrt{2C(\rho_t)}dW_t,\quad u_0\sim\pi_0,\]
where the macroscopic mean and covariance operator are defined as
\[m(\rho) = \int_{\mathcal{X}} x\rho(x)\, dx,\quad C(\rho) = \int_{\mathcal{X}} (x-m(\rho))\otimes(x-m(\rho))\, dx. \]
This connects the interacting Langevin system to its origin Langevin diffusion \eqref{eq:Langevin}. Hence, in the long-time limit the preconditioning matrix will formally be given by the covariance operator corresponding to the stationary distribution (assuming it exists).

The resulting modified Fokker--Planck equation is given by
\begin{equation}\label{eq:mod_FP}
\partial \rho_t = \nabla\cdot (\rho_t C(\rho_t) \nabla\log\pi_\ast) + {\mathrm{Tr}}(C(\rho_t)D^2\rho_t), \quad \rho_0 = \pi_0.
\end{equation}
Assuming that $C(\rho_t)\ge \alpha {\mathrm{Id}}$ and the target distribution of the form $\pi_\ast(u)\propto \exp(-\Phi(u))$, $\Phi(u) = \frac12\|y-\mathcal G(u)\|_\Gamma^2+\lambda \|u\|_{C_0}^2$, to be $\log$-concave, the solution $\rho_t$ of \eqref{eq:mod_FP} converges exponentially fast to equilibrium
$$
\mathrm{KL}(\rho_t\mid\pi_\ast)\le  \exp(-\lambda t) \mathrm{KL}(\rho_0\mid\pi_\ast),
$$
for some $\lambda>0$ \cite[Proposition~3.1]{GFLA20}. Furthermore, through the preconditioning with the sample covariance the resulting scheme remains invariant under affine transformations \cite{GNR2020}
%\textbf{TBD: include Gaussian special case (along the FP equation solutions remain Gaussian) and Wassertein convergence/ GD argumentation in \cite{GFLA20}?}

\subsection{Ensemble Kalman sampler}

One of the attractive features of the EnKF as well as of EKI is its derivative-free implementation.  The basis of the ensemble Kalman sampler (EKS) is to build a modified interacting Langevin sampler avoiding to compute derivatives.  Let $\pi_\ast(u)\propto \exp(-\frac12\|y-\mathcal G(u)\|_\Gamma^2-\|u\|_{C_0}^2)$, then the interacting Langevin system is given by 
\[du_t^{(j)} = -C(u_t){\mathrm D}\mathcal G(u_t^{(j)})^\top \Gamma^{-1}(G(u_t^{(j)})-y)-C(u_t)C_0^{-1}u_t^{(j)}\, dt + \sqrt{2C(u_t)}dW_t, \quad j=1,\dots,J.\]
Motivated by the approximation $C^{uw}(u) \approx C(u){\mathrm D}\mathcal G(u^{(j)})^\top$ the EKS is then formulated as the solution of the system of coupled SDEs
\[ du_t^{(j)} = -C^{uw}(u_t)\Gamma^{-1}(\mathcal G(u_t^{(j)})-y)-C(u_t)C_0^{-1}u_t^{(j)}\, dt + \sqrt{2C(u_t)}dW_t, \quad j=1,\dots,J.\]
We note that the approximation $C^{uw}(u) \approx C(u){\mathrm D}\mathcal G(u^{(j)})^\top$ is exact for linear forward models and hence, the EKS coincides with the interacting Langevin sampler in the linear setting. However, for nonlinear forward models the approximation of derivatives is only accurate in case the particles are close to each other. Since in the application of EKS the particles are aiming to represent a distribution, the particles are not expected to be close to each other. This fact suggests to formulate a localized version of the preconditioning sample covariance matrix, incorporating more weights on particles close to each other, but reducing the weight between particles far away.  Therefore, we define the distance-dependent weights between particle $u_t^{(j)}$ and $u_t^{(i)}$
\[w_t^{ji} = \frac{\exp(-\frac{1}{2\gamma}\|u_t^{(j)}-u_t^{(i)}\|_D^2)}{\sum_{l=1}^J\exp(-\frac{1}{2\gamma}\|u_t^{(j)}-u_t^{(l)}\|_D^2)},\]
for scaling parameters $\gamma>0$ and symmetric positive-definite matrix $D\in \mathbb{R}^{n_u\times n_u}$. The localized (mixed) sample covariance matrix around particle $u_t^{(j)}$ is the defined as
\[C(u_t^{(j)}) = \sum_{i=1}^J w_t^{ji} (u_t^{(i)}-\bar u_t^{(j)})\otimes (u_t^{(i)}-\bar u_t^{(j)}),\quad C^{uw}(u_t^{(j)}) = \sum_{i=1}^J w_t^{ji} (u_t^{(i)}-\bar u_t^{(j)})\otimes (\mathcal G(u_t^{(i)})-\bar{\mathcal G}_t^{(j)}), \]
with localized mean
\[ \bar u_t^{(j)} = \sum_{i=1}^J w_t^{ji} u_t^{(i)},\quad \bar{\mathcal G}_t^{(j)} = \sum_{i=1}^Jw_t^{ji} \mathcal G(u_t^{(i)}).\]
The localised EKS then reads as 
\[ du_t^{(j)} = -C^{uw}(u_t^{(j)})\Gamma^{-1}(\mathcal G(u_t^{(j)})-y)-C(u_t^{(j)})C_0^{-1}u_t^{(j)}\, dt + \sqrt{2C(u_t^{(j)})}dW_t, \quad j=1,\dots,J.\]
While the original EKS shows promising results for nearly Gaussian target distribution, the considered localized variant helps to extend the scope to multimodal target distributions \cite{RW2021}.
Other such related work has aimed to provide further understandings of the EKS. This has included the derivation of providing mean field limits and , but also providing various generalizations
\cite{DL21,GNR2020}.

\section{Other directions}
\label{sec:other}

As we have discussed some of the more recent developments in EKI, we now focus
on other, more smaller, extensions. In this section we will discuss these each in turn,
which will include machine learning, understanding EKI in the context of nonlinear
inverse problems, and finally applications related to engineering such as geophysical
sciences.

\subsection{Applications in machine learning}

The developments of machine learning methodologies has seen a significant increase in the last decade,
which have been produced to solve problems related to health-care, imaging, and decision processes. In particular
much of the these developments has been to due the advancements in optimizaion theory. As a result, ensemble
Kalman methods can be viewed as a natural class of algorithms to be directly applied, as they are derivative-free 
optimizers. 

The first work aimed at characterizing this connection was \cite{KS19} which demonstrated this. The authors motivated EKI
as a replacement to SGD where they initially applied it to supervising learning problems. Given a dataset $\{x_j,y_j\}_{j=1}^N$
assumed to be i.i.d. samples from a particular distribution, then given the Monte Carlo approximation one has the minimization
procedure
\begin{align}
\nonumber
&\arg \min_{u} \Phi_{\textrm{s}}(u;x,y), \\
\label{eq:lsf2}
\Phi_{\textrm{s}}(u;x,y) &= \frac{1}{N}\sum^N_{i=1}\mathcal{L}(\mathcal{G}(u|x_j),y_j) + \frac{\lambda}{2}\|u\|^2_{C_0},
\end{align}
where $\mathcal{L}:\mathcal{Y} \times \mathcal{Y} \rightarrow \mathbb{R}^+$, is some positive-definite function. 
In other words, one is trying to learn $x_j$ from the labelled data $y_j$. Supervised learning is used for common ML applications
such as image classification and natural language processing. Another related application is that of semi-supervised learning, which
aims to learn $x_j$ from some of the data $y_j$ where do not have access to all of it. This modified the least squares functional
given in \eqref{eq:lsf2}.

Another interesting direction has been the inclusion of EKI for training and learning neural networks \cite{HLR18}. This builds upon the previous
work discussed, but with a number of modifications. In particular what the authors show is that they are able to prove convergence of EKI
to the minimizer of a strongly convex function. They apply their modified methodology to a nonlinear regression problem of the form
$$
F(\theta) = A \theta + \epsilon \sin(B \theta),
$$
where $\theta$ is the parameter of interest and $F(\theta)$ is the objective functional of interest. This was also extended to the likes of
image classification problems, specifically the well-known MNIST handwritten data set.

A final and more recent direction of EKI and ML, was the work of Guth et al. \cite{GSW20}, which provided a way of solving the forward problem, 
within EKI. 

\subsection{Extensions to nonlinear convergence analysis}

A major challenge with EKI, and the EnKF in general, is establishing convergence analysis and properties in the nonlinear setting.
As it is well known in the linear and Gaussian setting, as the the number of particles $N \rightarrow \infty$, the EnKF
coincides with the KBF. However in the nonlinear setting it is has been challenging to derive any such results rigorously. 
Some ongoing and recent work has aimed to bridge the connections between EKI and nonlinear dynamics. The first paper that provided
some form of analysis was the work of Chada et al. \cite{CT21} which considered a specific form of EKI, in the discrete-time setting.

Namely the update formula is modified to
\begin{align*}
m_{n+1} &= m_n + C^{pp}_n(C^{up}_n + h_n^{-1} \Gamma)^{-1}(z-H(m_n)), \\
C_{n+1} &= C^{uu}_n - C^{up}_n(C^{pp}_n + h_n^{-1}\Gamma)^{-1}C^{pu}_n + \alpha^2_n \Sigma,
\end{align*}
where we adopt an ensemble square root filter formulation, which is known to perform better. As well as this we also
include covariance inflation (i.e. inflation factor of $\alpha_n$), and an adaptive step-size $h_n$ motivated from stochastic optimization to allow an acceleration
for the convergence. However the other underlying contribution, as eluded to, is that given this update form we are able
to prove convergence towards both local and global minimizers. In other words for the later, we have the following result
$$
\lambda_c\|m_N - u^*\|^2 \leq \ell(m_N) - \ell(u^*) \leq \frac{D}{N^{\alpha}},
$$
which the above result establishes polynomial convergence. We note from the above equation, that $\lambda_c$ is a convexity constant,
$\ell$ is the associated loss function, $D$ is some constant, $u^*$ is the global minimizer and $\alpha$ is some term, which we refer
to \cite{CT21}, for further details.

As one can notice, this convergence analysis was considered for the discrete-time setting, so a natural extension from this is to the continuous-time
framework. The work of Blomker et al. \cite{BSW21} provide a first convergence analysis in this direction. However given both these works, a full understanding
in the nonlinear setting has not been achieved, where considerable work is still required. Thus these papers provide a first step in doing so, for both settings.

\subsection{Engineering applications}

As a final direction to discuss in detail, which is very much related to the theme of this book, are applications in particular engineering applications.
The advantage of these ensemble Kalman methods, is that they can be viewed as a black box-solver, therefore it is highly applicable. One particular
application has been geophysical sciences, related to recovering quantities of interest which are below the surface, or subsurface. Examples include
the inverse problem of electrical resistivity tomography (ERT), shown below 
\begin{figure}[!h]\centering
	\FIG{\includegraphics[width=0.5\textwidth,height=9pc]{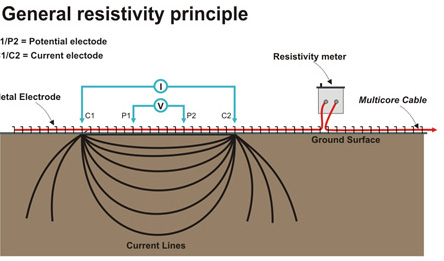}}
	{\caption{Image depicting electrical resistivity tomography, where the the electric currents are recorded at the electrodes of the subsurface material.}}
	\label{fig:ERT}
\end{figure}

ERT is concerned with recovering, or characterizing sub-surface materials in terms of their electrical properties, which are recorded through electrodes.
It operates very similarly to electrical impedance tomography (EIT), expect the difference being that it is subsurface. This has been also considered for
learning permeability of subsurface flow in a range of different settings which can be found in the following papers \cite{TIW21,MT20}. 
\\
Another interesting direction is related to walls, specifically quantifying uncertainty in thermo-physical properties of walls. This work was conducted by Iglesias et al
\cite{DIJ18,ISS18}. Specifically the application is the inverse problem of recovering the thermodynamic property or temperature. Similar work related to the methodology used
here has been used in resin transfer modeling \cite{IPT18}, based on problems of moving boundaries. This is a difficult problem to model, however it provides a 
first step in doing so. Aside from these applications other particular applications include mineral exploration scattering problems, 
numerical climate models and others. \cite{SAS21,DGS21,HLW21}. It is worth mentioning that, as of now, there is no official online software package for EKI in general.
This is currently being developed, but we emphasize to the reader that the methodology presented, with the examples later, are not related to well known softwares that
are available in \texttt{Matlab} or \texttt{Python}.

As a side remark, there are more directions beyond what is discussed above. Some others, without going into details, include developing hierarchical approaches, 
 incorporating constrained optimization, and connections with data assimilation strategies \cite{AB19, CCS21, CIRS18, CSW19, TM22}.

\section{Numerical experiments}
\label{sec:num}

In this section we provide some numerical experiments highlighting the performance of ensemble Kalman methods for inverse problems.
Specifically we will consider EKI as discussed in Section \ref{sec:EKI_back}.  We will compare EKI with its regularized version of TEKI. Both
these methodologies will be tested on on two motivating inverse problems arising in geophysical and atmospheric sciences, i.e. a Darcy flow
partial differential equation and the Navier--Stokes Equation.

In order to assess a comparison, we will present three different figures. (i) The first being a reconstruction at the end 
of the iterative scheme; (ii) the error between the approximate solution and the ground truth, and (iii) the data misfit. The equations associated with each are given as
\begin{itemize}
\item \textit{Reconstruction through EKI}: $\frac{1}{J}\sum^J_{j=1}u^{(j)}_n$.
\item \textit{Relative error}: $\frac{\|u^{\dagger} - u\|^2_{L^2}}{\|u^{\dagger}\|_{L^2}}$.
\item \textit{Data misfit}: $\| \Gamma^{-1/2}(y-\mathcal{G}(u^{\dagger})\|^2$.
\end{itemize}

\subsection{Darcy flow}

Our first model problem is an elliptic partial differential equation (PDE), which has numerous applications. Specifically one of them is subsurface flow
in a porous medium. The forward problem is concerned with solving for the pressure $p \in H^1_0(\Omega)$, given the permeability $\kappa \in L^{\infty}(\Omega)$
and source function $f \in L^{\infty}(\Omega)$, where the PDE is given as
\begin{align}
\label{eq:pde}
-\nabla \cdot (\kappa \nabla  p) &= f, \quad \in \Omega, \\
\label{eq:bc}
p&=0, \quad \mathrm{on} \ \Omega.
\end{align}
such that we have prescribed Dirichlet boundary conditions, and $\Omega = [0,1]^2 \subset \mathbb{R}^d$, for $d=2$, is a Lipschitz domain. The inverse problem associated
to solving $p$ from \eqref{eq:pde} is the recovery of the permeability $\kappa \in L^{\infty}(\Omega)$, from noisy measurements of $p$, i.e.
$$
y = \mathcal{G}(\kappa) + \eta, \quad \eta \sim \mathcal{N}(0,\Gamma),
$$
where recalling that $\mathcal{G}(\kappa) = p$. We consider 64 equidistance observations within the domain, and on the boundary. 
To numerically solve \eqref{eq:pde} we employ a centered-finite difference method with a mesh size of $h=1/100$. For our
noisy observations we consider $\Gamma = \gamma I$, where $\gamma=0.01$. We will use and compare EKI and TEKI, with an ensemble size of $J=50$ for both methods.
We will run both iterative schemes for $n=24$ iterations.  For our initial ensemble $\{u_0\}^{J}_{j=1}$ we consider modelling it as a Gaussian random field, i.e. $u \sim N(0,C)$, which can be done
via the Karhunen-Lo\`{e}ve expansion
\begin{equation}
\label{eq:kle}
u = \sum_{k \in \mathbb{Z}^+} \sqrt{\lambda_k}\phi_k \xi_k, \quad \xi_k \sim N(0,1),
\end{equation}
where $(\lambda_k,\phi_k)$ are the associated eigenvalues and eigenvectors of the covariance operator $C$. There are numerous choice of covariance functions one can take, however a popular
choice is the Mat\'{e}rn covariance function, which provides much flexibility for modelling. For full details on various covariance functions, or operators, we refer to reader to \cite{LPS14}.
The true unknown of interest is taken to be also a Gaussian random field, but one that is smoother than that of that of the initial ensemble. 

Our first set of experiments are provided in \textbf{Figure 4} which shows the truth, the reconstruction from using EKI, and that of using TEKI. As we can observe, is it clear that both methodologies
work well at learning the true unknown function. However it is clear that the TEKI induces a smoother reconstruction, which arises from the regularization. However, what is interesting is that
if we analyze \textbf{Figure 5}, we notice that the relative error tends to diverge at the end with EKI, and this is due to the overfitting of data. A motivation behind TEKI is to alleviate this. This can be seen
vividly as it tends to decrease, and for the data misfit, it remains within the noise level, which is given as
$$
\textrm{noise \ level} \  = \|(y - G(u^{\dagger}))\| = \| \eta^{\dagger}\|.
$$

\begin{figure}[!h]\centering
	\FIG{\includegraphics[width=0.9\textwidth]{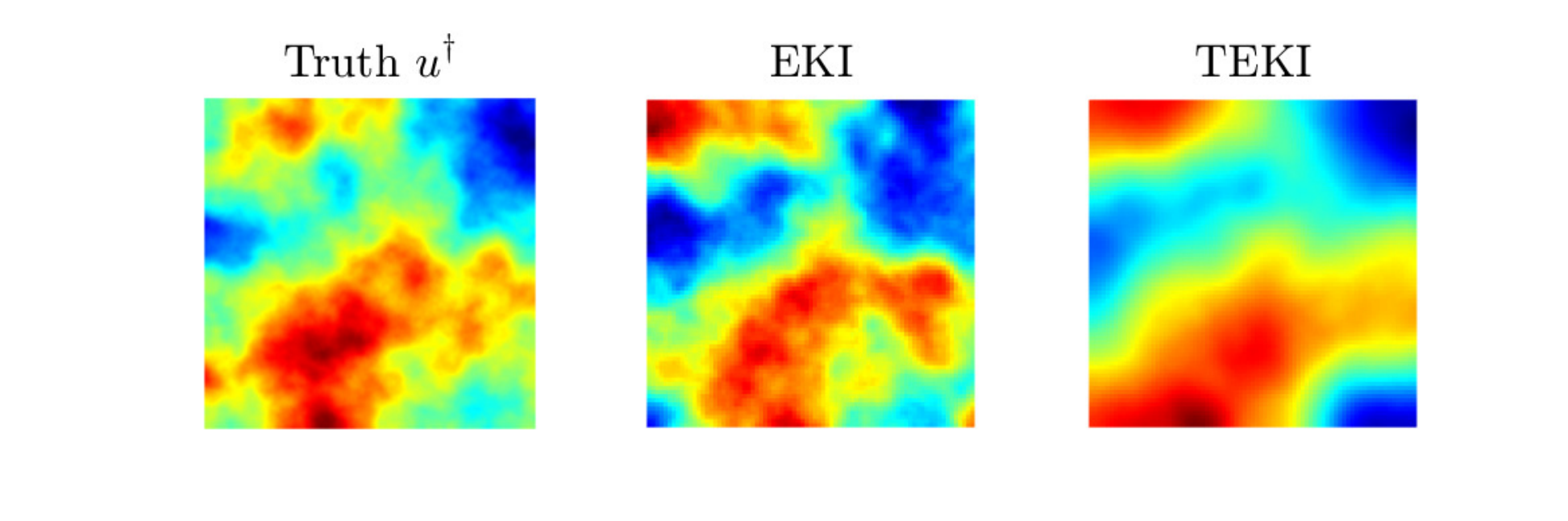}}
	{\caption{Reconstruction plots for the Darcy flow PDE example. Left: Truth. Middle: EKI reconstruction. Right: TEKI reconstruction.}}
	\label{fig:pde_recon}
\end{figure}

\begin{figure}[h!]\centering
\includegraphics[scale=0.31,trim=0 0 0 0]{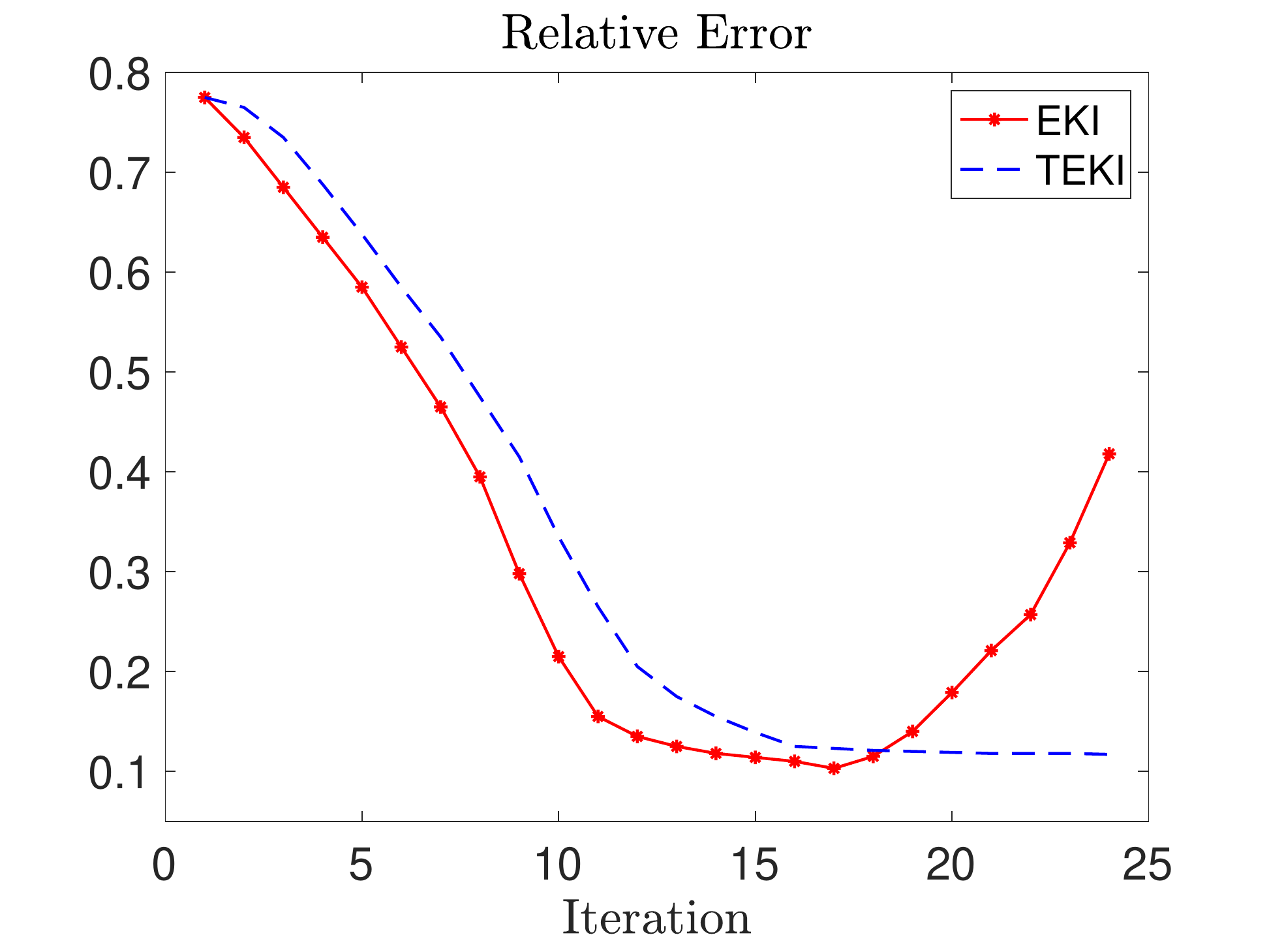}
\vspace{3mm}
\includegraphics[scale=0.31,trim=0 0 0 0]{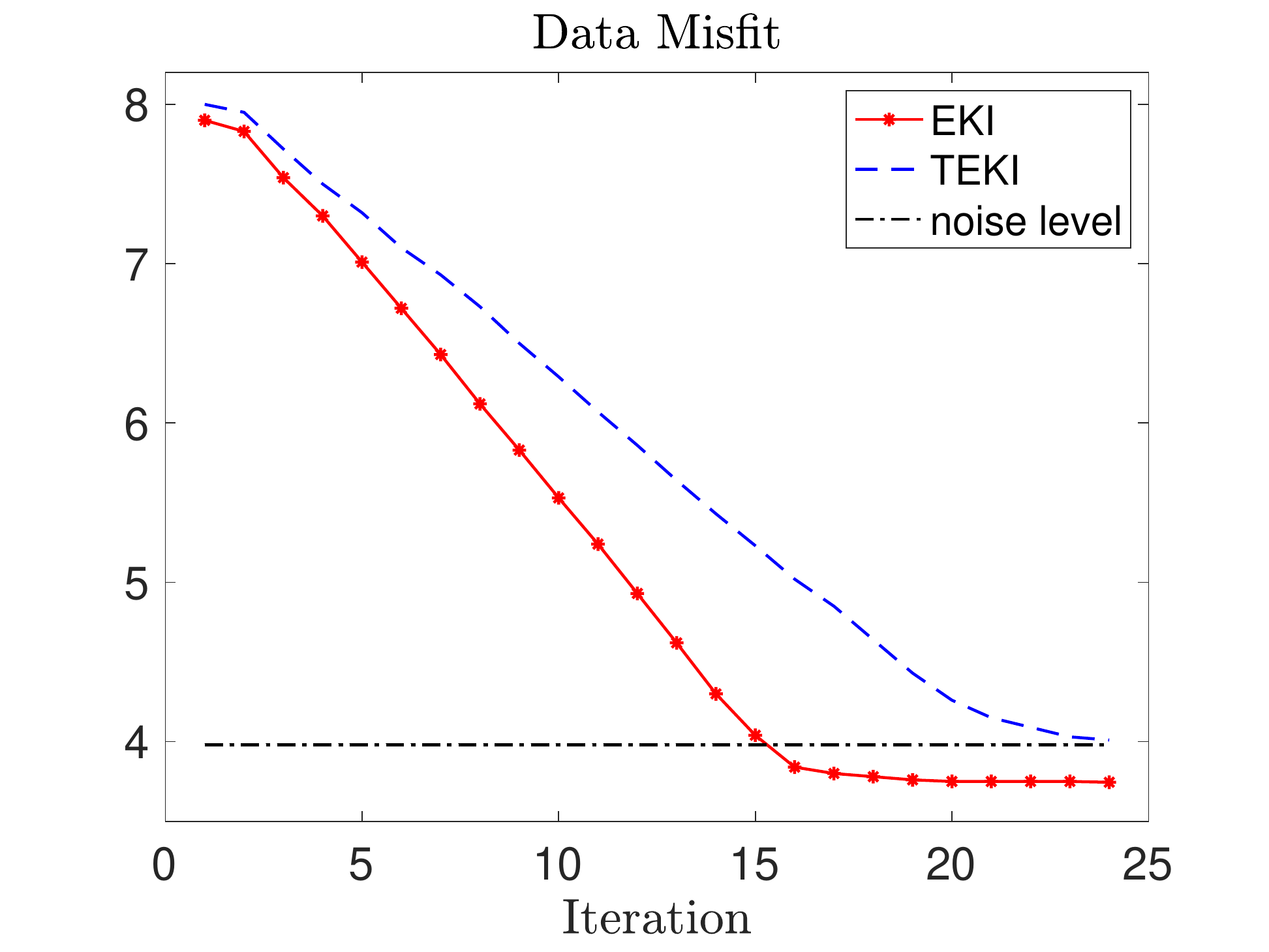}
	{\caption{Relative errors and data misfits for the Darcy flow PDE example. We compare EKI with TEKI.}}
	\label{fig:pde_error}
\end{figure}

\subsection{Navier--Stokes equation}

\begin{figure}[!h]\centering
	\FIG{\includegraphics[width=0.9\textwidth]{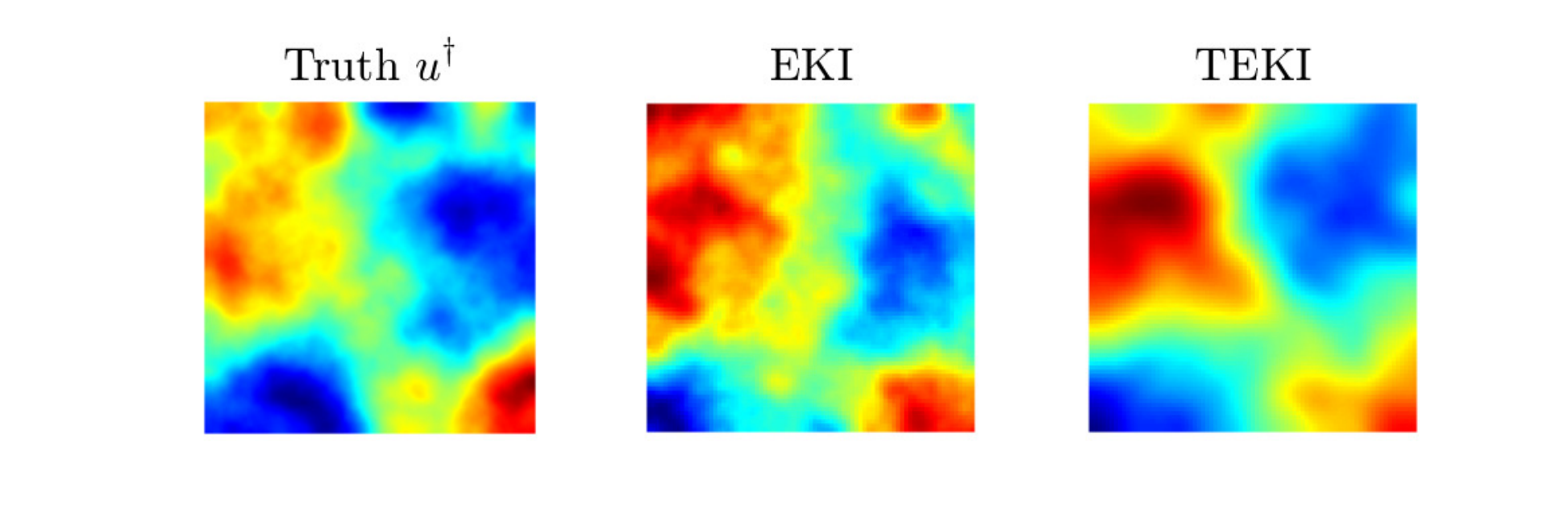}}
	{\caption{Reconstruction plots for the NSE PDE example. Left: Truth. Middle: EKI reconstruction. Right: TEKI reconstruction.}}
	\label{fig:nse_recon}
\end{figure}

Our final test problem is a well-known PDE model arising in numerical weather prediction
which is the Navier--Stokes equation (NSE).
We consider a 2D NSE defined on a torus $\mathbb{T}^2 = [0,1]^2$ with periodic boundary conditions. The aim to estimate the velocity $v:=[0,\infty) \times \mathbb{T}^2 \rightarrow \mathbb{R}^2$ defined as a vector field from the scalar pressure field $p:=[0,\infty) \times \mathbb{T}^2 \rightarrow \mathbb{R}^2$. The NSE is given as 

\begin{align}
\label{eq:NSE}
\partial_tv+ (v \cdot \nabla)v + \nabla p - \nu \Delta v &= f,  \quad [0,\infty) \times \mathbb{T}^2,  \\
\label{eq:flux}
\nabla \cdot v &=0, \quad [0,\infty) \times \mathbb{T}^2, \\
\label{eq:IC}
v &=u, \quad \{0\} \times \mathbb{T}^2, 
\end{align}
with initial condition \eqref{eq:IC} and zero flux \eqref{eq:flux}. From \eqref{eq:NSE} $f\in [0,\infty) \times \mathbb{T}$ corresponds to a volume forcing, 
$\nu$ is the associated viscosity of the fluid. For the NSE equation we consider a spectral Fourier solver for \eqref{eq:NSE}. The PDE is more challenging 
to invert than the previous example, therefore we take 100 point-wise observations. The setup is largely the same as the previous example, where we take
an initial condition based on a Gaussian random field through the KL expansion \eqref{eq:kle}. We will aim to recover the true underlying function $u^{\dagger}$
using both EKI and TEKI.  The results are obtained from the experiments are presented in \textbf{Figure \ref{fig:nse_recon}}  and \textbf{Figure 7}. A similar phenomenon shows,
where the reconstructions work well, however there is an additional smoothness induced through the regularization in TEKI. Similarly, as we see with the relative errors
and data misfit the overfitting of the data in the end for EKI. We note that this can be avoided depending on the prior form, its hyperparameters, the observations, and the noise. 
However we specify particular choices to demonstrate it can occur.

\begin{figure}[h!]
\centering
\includegraphics[scale=0.31,trim=0 0 0 0]{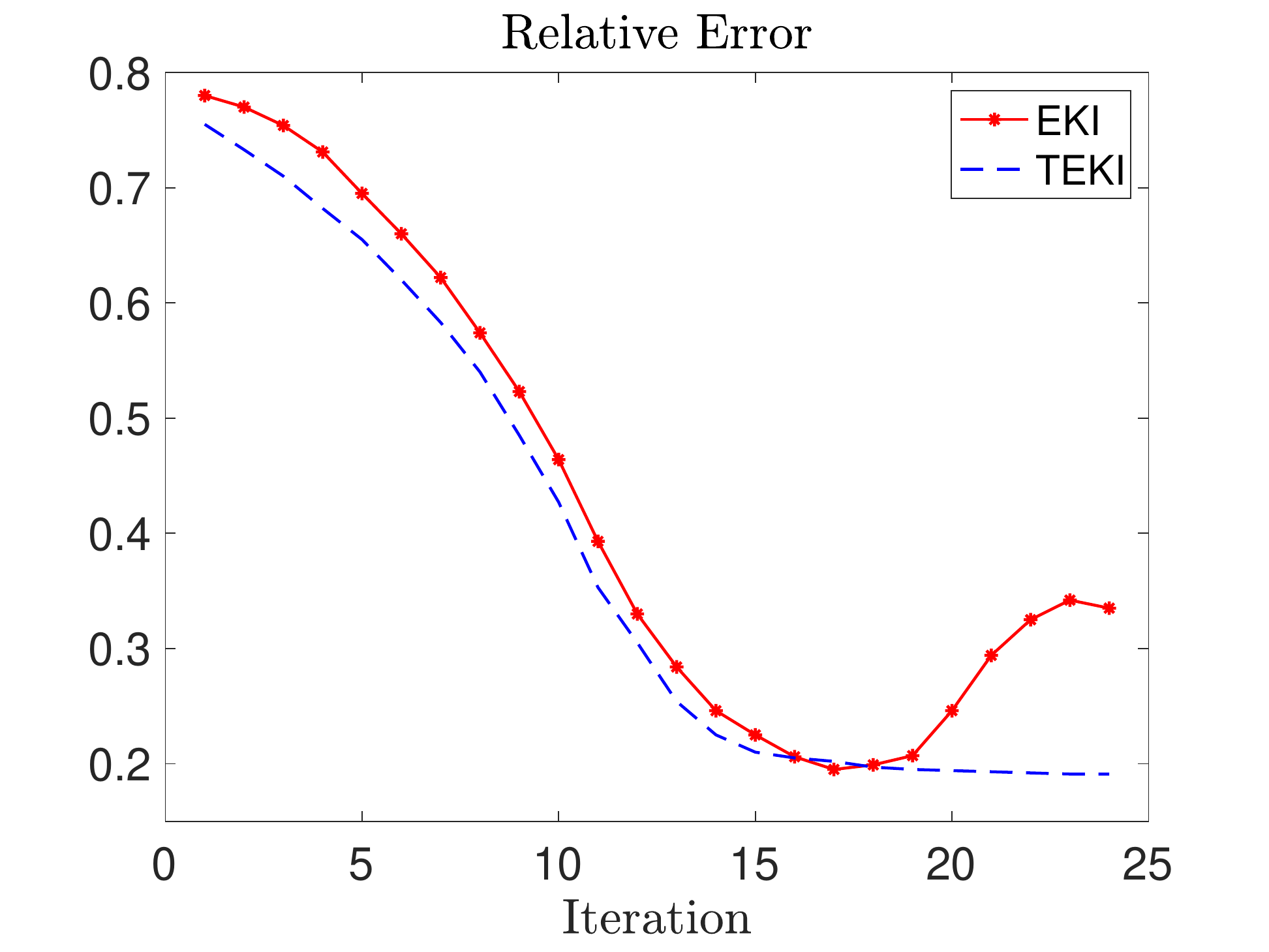}
\vspace{3mm}
\includegraphics[scale=0.31,trim=0 0 0 0]{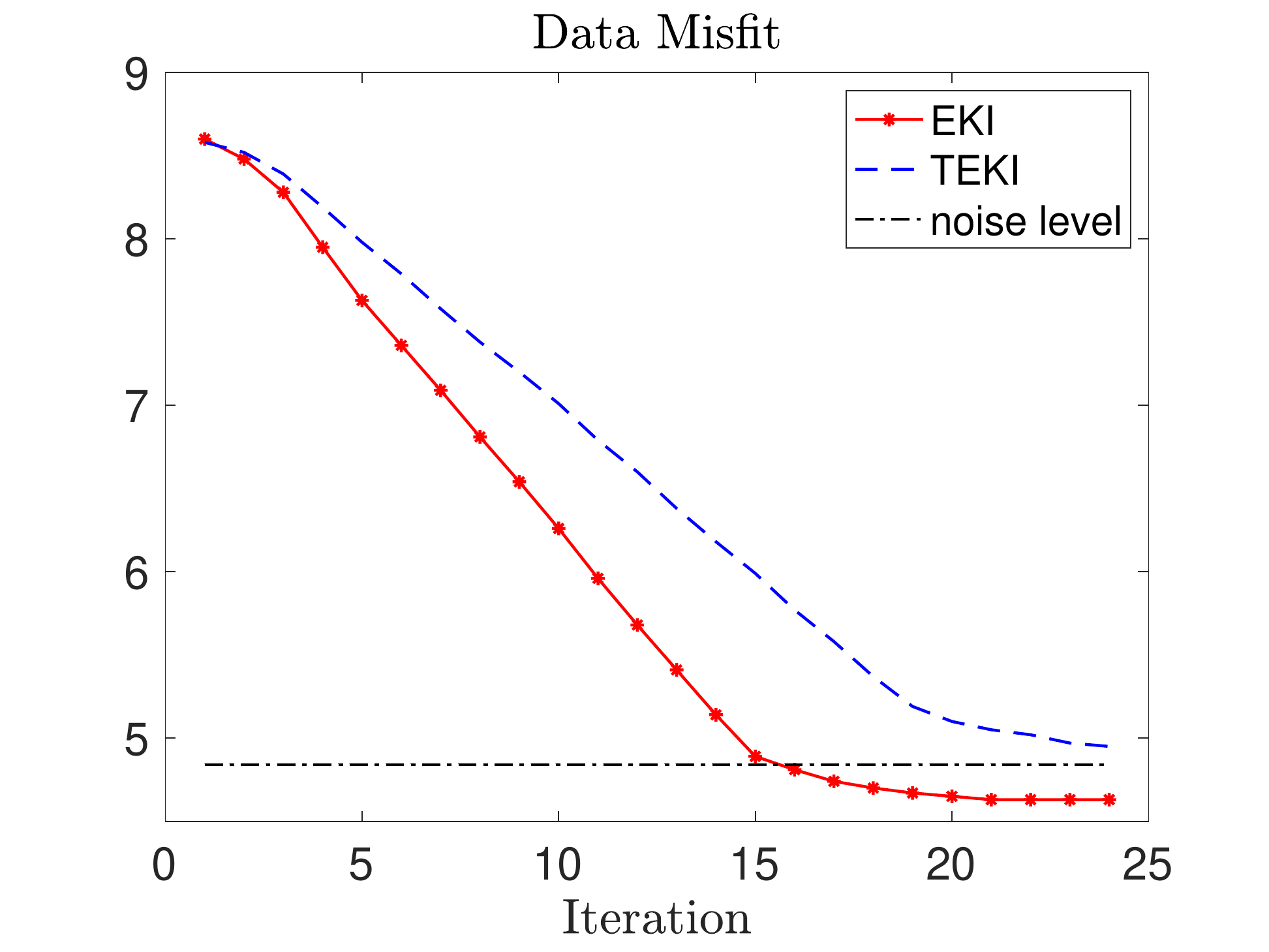}
	{\caption{Relative errors and data misfits for the NSE PDE example. We compare EKI with TEKI.}}
	\label{fig:nse_error}
\end{figure}

\section{Conclusion}
\label{sec:conc}

The ensemble Kalman filter (EnKF) is a simplistic, easy-to-implement and powerful algorithm.
This has been particularly the case in numerous data assimilation applications for state estimation,
which includes the likes of numerical weather prediction, geosciences and more recently machine learning. 
A major advantage of the method is that, unlike other filters such as the particle filter, it scales better in high
dimensions, and can be significantly cheaper. In this chapter we consider the EnKF and its application to
parameter estimation. Such a mathematical procedure also has similar applications to the ones states, where
one can exploit such techniques for inverse problems. We provide a review and overview of some of the major
contributions in this direction, where the resulting methodology is known as ensemble Kalman inversion (EKI),
based largely on the work of Iglesias et al. \cite{ILS13}.
We presented various avenues the field of EKI has taken such as regularization, extensions to sampling, and other
areas. We demonstrated how EKI can perform on two numerical examples PDE examples.

The EKI methodology is one which builds very naturally from many different fields, which acts a strong motivation. For example
being an optimizer, one can naturally apply optimization procedures, but also techniques from data assimilation and uncertainty
quantification. As a result, this methodology naturally brings researchers from different fields working towards parameter estimation,
and inverse problems. This synergy of areas will hopefully ensure new emerging directions within EKI, from a methodological,
theoretical and application perspective.

\section*{Acknowledgments}
This work was funded by KAUST baseline funding. The author thanks Simon Weissmann for helpful discussions, and for the use of some of the earlier figures,
and information on EKS.

\section*{Abbreviations}

\begin{abbrvlist}[DMEM-FBS]
	\item[EnKF] Ensemble Kalman filter
	\item[EKI] Ensemble Kalman inversion
	\item[EKS] Ensemble Kalman sampler
	\item[EIT] Electrical impedance tomography
	\item[ERT] Electrical restivitiy tomography
	\item[MNIST] Modified National Institute of Standards and Technology 
	\item[KL] Kullback-Leibler
	\item[KF] Kalman filter
	\item[L96] Lorenz 96 model
	\item[LSF] Least squares functional
	\item[PDE] Partial differential equation
	\item[SDE] Stochastic differential equation
\end{abbrvlist}

%\section*{Appendix A}
%The philosophy behind the score sheets states that individuals or groups define
%a standard for the quality of their performance. Then, they describe the standard in
%terms of a set of requirements. This set is the score sheet. It allows for peer evaluation
%and self-evaluation of an activity. Grading proceeds by determining the
%fraction of requirements fulfilled and is objective and reproducible. The score sheet
%exists prior to the execution of any activity and thus induces iteration until the
%performance becomes satisfactory. Any individual or group can adapt the method
%to any professional activity by selecting the pertinent requirements to fulfill their
%standard of excellence.

\begin{authordetails}
	
	% Author details will always appear the end of the chapter in the final version of the chapter
	
	\author{Neil K. Chada$^{1}$}
	%\author{Simon Weissmann$^{2}$}
	%
	\address[1]{King Abdullah University of Science and Technology, Thuwal, 23955, Saudi Arabia}
	%\address[2]{Interdisciplinary Center for Scientific Computing, University of Heidelberg, 69120, Heidelberg, Germany}
	%\address[3]{Institution No. 3, City, Country}
	%
	\address{*Address all correspondence to: neilchada123@gmail.com}
	%
	%\address{\dag\ These authors contributed equally}
	
	\IntechOpentext{\textcopyright\ \the\year{} The Author(s). License IntechOpen. This chapter is distributed under the terms of the Creative Commons Attribution License (http://creativecommons. org/licenses/by/3.0), which permits unrestricted use, distribution, and reproduction in any medium, provided the original work is properly cited.}
	
	% Note: The copyright year will be changed accordingly during production to correspond with the year of publication.
	
\end{authordetails}

\end{document}